\documentclass[leqno,10pt]{amsart}
\usepackage{graphicx}
\usepackage{float}

\usepackage{amsmath,amsthm,amsbsy,amsfonts,amssymb}
\usepackage[backref=page]{hyperref}

\hypersetup{colorlinks=true,citecolor=blue,linkcolor=blue,urlcolor=blue,pdfstartview=FitH,
pdfauthor=Beno\^it Louvel, pdftitle=The first moment of Sali\'e sums}

\theoremstyle{plain}
\numberwithin{equation}{section}
\newtheorem{thm}{Theorem}[section]
\newtheorem{lemme}[thm]{Lemma}
\newtheorem{prop}[thm]{Proposition}
\newtheorem{cor}[thm]{Corollary}

\theoremstyle{definition}
\newtheorem*{ack}{Acknowledgements}

\newtheorem{defi}[thm]{Definition}

\theoremstyle{remark}
\newtheorem{rque}[thm]{Remark}

\begin{document}

\newcommand{\rond}{\mathcal}
\newcommand{\id}{ \mathfrak}
\newcommand{\einb}{\hookrightarrow} 
\newcommand{\too}{\longrightarrow} 
\newcommand{\pt}{\cdot} 
\newcommand{\pts}{\ldots}
\newcommand{\N}{\mathbb{N}}
\newcommand{\Z}{\mathbb{Z}}
\newcommand{\Q}{\mathbb{Q}}
\newcommand{\R}{\mathbb{R}}
\newcommand{\A}{\mathbb{A}}
\newcommand{\F}{\mathbb{F}}
\newcommand{\C}{\mathbb{C}}
\newcommand{\gras}{\textbf}
\newcommand{\dual}{\land}
\newcommand{\ec}{\textnormal}
\newcommand{\tq}{\vert}
\newcommand{\iso}{\cong}
\newcommand{\tild}{\widetilde}
\newcommand{\congru}{\equiv}
\newcommand{\notcongru}{\not\equiv}
\newcommand{\barre}{\overline }
\newcommand{\norme}{\Vert}
\newcommand{\Mmod}{\ (\bmod\ }
\newcommand{\mmod}{\ (\!\bmod }

\title{The first moment of Sali\'e sums}

\author{Beno\^it Louvel}
\address{Mathematisches Institut G\"ottingen \\
Bunsenstr. 3-5 \\ 37073 G\"ottingen}
\email{blouvel@uni-math.gwdg.de}

\thanks{Research supported in part by the Volkswagen Fundation}
\begin{abstract} 
The main objective of this article is to study the asymptotic behavior of Sali\'e sums over arithmetic progressions. We deduce from our asymptotic formula that Sali\'e sums possess a bias of being positive. The method we use is based on Kuznetsov formula for modular forms of half integral weight. Moreover, in order to develop an explicit formula, we are led to determine an explicit orthogonal basis of the space of modular forms of half-integral weight. 
\end{abstract}
\maketitle
\setcounter{page}{1}

\section{Introduction}\label{introduction}

For an odd integer $c$, the Sali\'e sum is defined as 

\begin{equation*}
K_2(m,n;c)=\sum_{\substack{x\Mmod c) \\ x x^{-1} \congru 1 \Mmod c)}} \left(\frac{x}{c}\right) e\left(\frac{mx+nx^{-1}}{c}\right).
\end{equation*}
As usual, for $z\in\C$, we write $e(z)=\exp(2\pi i z)$. We use the subscript $K_2$ in order to distinguish the Sali\'e sums from the Kloosterman sums, defined analogously, but without the presence of the Jacobi symbol $(x/c)$.\\ 

Classically, these sums have been investigated over the set of primes. For a prime $p$, one has, by a classical theorem of Sali\'e, that

\begin{equation}\label{intro:eq:salie}
K_2(m,n;p)= 2 \cos\left(\frac{4\pi x}{p}\right) \sum_{y\Mmod p)} e\left(\frac{y^2}{p}\right), \quad x^2 \congru mn \Mmod p). 
\end{equation}
From the equidistribution of roots of quadratic congruences modulo prime moduli (proved in \cite{dfi7})  and from \eqref{intro:eq:salie}, follows that the angles of Sali\'e sums $K_2(m,n;p)$ are equidistributed with respect to the uniform measure.\\ 

In this paper, we investigate the distribution of the sums $K_2(m,n;c)$, when $c$ runs over the set of integers.
Individually, these sums are well understood. One can generalize \eqref{intro:eq:salie} to any modulus $c$ and show that Sali\'e sums satisfy the individual bound

\begin{equation}\label{intro:eq:bound}
\left\vert K_2(m,n;c)\right\vert \leqslant 2^{\omega(c)} \sqrt{c},
\end{equation}
where $\omega(c)$ is the number of distinct prime divisors of $c$. Nevertheless, as for most of (complete) exponential sums, the understanding of the behavior of these sums as function of the moduli $c$ is a difficult problem. Using methods of analytic number theory, one might expect to catch important properties of Sali\'e sums by looking at their $L$-function, but it turns out that, due to a twisted multiplicativity, the $L$-function has neither a functional equation nor an Euler product, rendering the investigation of Sali\'e sums very involved.\\  

One way these sums can be approached is by the analytic theory of automorphic forms. It is known that Kloosterman sums (or their twists) appear as Fourier coefficients of Poincar\'e series. In 1987, in his breakthrough paper \cite{iwa:fourier-coeff}, Iwaniec succeded in proving a new upper bound for the Fourier coefficients of modular forms of half-interal weight, by estimating sums on Sali\'e sums. This leads naturally to the problem of studying Sali\'e sums on average, and more precisely to detect cancellation among these sums. Inspired by the work of Livn\'e and Patterson \cite{liv-pat}, where the authors study the first moment of cubic exponential sums, we obtain a complete determination of the first moment of Sali\'e sums over arithmetic progressions. A crucial point for us is that the $L$-function associated to Sali\'e sums possesses an {\it exceptional pole}, related to the minimal eigenvalue of the Laplacian operator. It should be emphasised that this phenomenon is reminiscent to the situation in \cite{liv-pat}, and ultimately lies at the heart of the problem. \\

We study the distribution of the normalized Sali\'e sums $K_2(m,n;c)$, for fixed $m$ and $n$, when $c$ runs over an arithmetic sequence $c\congru 0 \Mmod D)$. One of the principal consequences of our main theorem is that, in most cases, the Sali\'e sums $K_2(m,n;c)$ exhibit much cancellation and that, in the remaining cases, our formula shows a definitive bias for the Sali\'e sums to being positive. . \\
 
As usual, we write $d \mid b^\infty$ if $d$ is supported by $b$, i.e. if $d$ is a product of primes dividing $b$, and $(a,b^\infty)$ for the greatest divisor of $a$ supported by $b$. For a positive integer $d$, we denote by $\chi_d$ the primitive quadratic character corresponding to the field extension $\Q(\sqrt{d})/\Q$. 
Let us define the symbol $\epsilon_c$, for odd integers $c$, by $\epsilon_c=1$ if $c\congru 1 \Mmod 4)$, and $\epsilon_c=i$ if $c\congru 3 \Mmod 4)$. One easily sees that the numbers $K_2(m,n;c) \barre{\epsilon_c}$ are real.
 
\begin{thm}\label{intro:thm}
Let $D$ and $f$ be odd positive integers, mutually coprime and let $\chi$ be an even primitive Dirichlet character modulo $f$. Let $m,n\in \Z$. Let $\varepsilon>0$ and let $X\gg 1$. 

If $\chi\neq \chi_f$, or if $m<0$ or if $n<0$, then

\begin{equation*}
\sum_{\substack{0<c\leqslant X \\ c\congru 0\Mmod D)}}
\frac{K_2(m,n;c)}{\sqrt{c}} \barre{\epsilon_c} \chi(c) = \rond{O}\left(X^{3/4+\varepsilon}\right).
\end{equation*}

If $\chi=\chi_f$ and $m,n>0$, then
 
\begin{equation*}
\sum_{\substack{0<c\leqslant X \\ c\congru 0\Mmod D)}}
\frac{K_2(m,n;c)}{\sqrt{c}} \barre{\epsilon_c} \chi(c) = C(D,f,m,n) X +\rond{O}\left(X^{3/4+\varepsilon}\right),
\end{equation*}
for some real number $C(D,f,m,n)$. Assume that $m$ and $n$ are of the form 

\begin{equation*}
\left\{\begin{aligned}
m&=tfs^2m'^2\\
n&=tfs^2n'^2,
\end{aligned}\right.
\end{equation*}
for some positive integers $t$, $s$, $m'$, $n'$ satisfying the following condions:

\begin{itemize}
\item[$(i)$] $t$ is square-free, $t\congru 1 \Mmod 4)$, $s$ is supported by $t$ and $s^2t^3 \tq D$,
\item[$(ii)$] $(m',t)=(n',t)=1$,
\item[$(iii)$] $(m',D_1)=(n',D_1)$, where $D=D_tD_0D_1^2$, where $D_t=(D,t^\infty)$ and $D_0$ is square-free.
\end{itemize}
Then, one has

\begin{equation*}
\begin{split}
C(D,f,m,n) = \frac{8}{\pi^2} \frac{s\sqrt{t}}{D} \left(\frac{fm'n'}{t}\right) (m',D_1) \prod_{p\tq Df} (1+p^{-1})^{-1} \\
\times \prod_{p\tq t}(1-p^{-1})^{-1} \prod_{p\tq D_1/(m',D_1)}(1-p^{-1})^{-1}.
\end{split}
\end{equation*}
Moreover, if $m$ and $n$ are not of the form described above, then $C(D,f,m,n)=0$. In particular, $C(1,f,m,n)\ge 0$.
\end{thm}

Let us now mention some other problems related to Theorem\,\ref{intro:thm}. Consider the sum

\begin{equation}
S(m,\ell;c)= \sum_{\substack{x\Mmod c)\\ x^2\congru m \Mmod c)}} e\left(\frac{2x\ell}{c}\right).
\end{equation}
These sums are related to Sali\'e sums by the formula

\begin{equation}\label{intro:eq:short-sum}
S(m,\ell;c) = \sum_{d\tq (\ell,c)} \sqrt{\frac{d}{c}} \, \barre{\epsilon}_{c/d} K_2\left(m,\frac{\ell^2}{d^2};c\right).
\end{equation}
From \eqref{intro:eq:short-sum} and Theorem\,\ref{intro:thm}, we obtain the following the Corollary.

\begin{cor}\label{intro:cor}
Let $m,\ell,c\in\Z$. Then for any $\varepsilon>0$, one has

\begin{equation*}
\sum_{c\leqslant X} S(m,\ell;c) = C(m,\ell) X +\rond{O}_{m,\ell}\left(X^{3/4+\varepsilon}\right),
\end{equation*}
with 

\begin{equation*}
C(m,\ell)=
\begin{cases}
0&\ec{if $m$ is not a square}\\
\frac{8}{\pi^2} X\frac{\sigma(\ell)}{\ell} &\ec{if $m$ is a square}
\end{cases}
\end{equation*}
Here, $\sigma(\ell)$ is the sum of the divisors of $\ell$.
\end{cor}
A similar result to Corollary\,\ref{intro:cor} is given in \cite[Theorem\,1]{hoo:div-quadr-poly}, where the author obtains a control on the dependance on $\ell$ of the error term, but only in the case where $m$ is not a square. 
Note that recently, Duke, Imamoglu and T\'oth have conjectured in \cite{dit:cycle-mock} that 

\begin{equation}\label{conj}
d^{-1/2} \sum_{c>0} S(d,m;c) \sin\left(\frac{4\pi m\sqrt{d}}{c}\right) \ll \sigma_1(m) Tr_d(1),
\end{equation}
where $Tr_d(1)$ is a trace of singular moduli; for example, if $d>1$ is a fundamental discriminant, then $Tr_d (1) = L(1,\chi_D)$.Theorem\,\ref{intro:thm} could also be applied to the asymptotic distribution of Dedekind sums, using the connection between Sali\'e sums and Dedekind sums; this has been studied for example by Vardi in \cite{var:these}. \\ 

Let us now give some indications on our proof of Theorem\,\ref{intro:thm}. For a congruence subgroup $\Gamma\subset \Gamma_0(4)$, one can define exponential sums $K_{\sigma,\tau}(m,n;c)$, which are associated to two cusps of $\Gamma$. The {\it geometric Sali\'e sums} $K_{\sigma,\tau}(m,n;c)$ appear as Fourier coefficients of non-holomorphic Poincar\'e series, and therefore one can relate sums of geometric Sali\'e sums with the spectrum of Maa\ss{} forms of weight $1/2$, using  the Kuznetsov trace formula. The geometric side of the trace formula is the easiest to deal with: we can choose in an appropriate way the cusps of the congruence subgroup $\Gamma=\Gamma_0(4Df)$, so that one can relate the geometric Sali\'e sums to the classical Sali\'e sums $K_2(m,n;c)$. The spectral side of the trace formula is more subtle to handle. A specific feature of Maa\ss{} forms of weight $1/2$ is to have an exceptional eigenvalue located at $\lambda=3/16$. We will show that square-integrable Maa\ss{} forms with respect to the exceptional eigenvalue $\lambda=3/16$ are in bijection with modular forms of weight $1/2$. 
As a result of independent interest (Theorem\,\ref{mod:thm}), we are able to obtain an explicit orthogonal basis for the space of modular forms of half-integral weight. This result then allows us to finally obtain an explicit expression for the spectral side of the trace formula. 

\begin{rque}
In principle, Theorem\,\ref{mod:thm} could be used to complete a result of Blomer in \cite{blo:hecke-quadr-poly}. There, the author studies the eigenvalues of Hecke eigenforms over quadratic polynomials and obtains a formula for their asymptotic behavior. In the case where the eigenform is chosen to be a Poincar\'e series, one can get an explicit expression for the asymptotic constant appearing in \cite[(1.4)]{blo:hecke-quadr-poly}, by means of our Theorem\,\ref{mod:thm}.
\end{rque}  

Finally, in order to illustrate the rate of convergence of sums of Sali\'e sums, some numerical examples for Theorem\,\ref{intro:thm} are presented in Section\,\ref{sec:num}. For a related discussion on this subject, concerning the distribution of Sali\'e sums and other arithmetical functions, we refer to \cite{pat:distrib-arith}.\\

\noindent{\it Notations} \
For a complex number $z$, we define its argument to be in the interval $[0,2\pi[$. The greatest common divisor of $a$ and $b$ is denoted by $(a,b)=\ec{gcd}(a,b)$, $\mu$ is the M\"obius function and $\varphi$ the Euler function. 
Let $\chi$ be a Dirichlet character of modulus $f$. By $\chi\chi_d$, we mean the primitive character associated to the product of $\chi$ and $\chi_d$. We denote the conductor of the primitive character associated to $\chi\chi_d$ by $f_d$; thus in particular, $f_1\mid f$, with equality if and only if $\chi$ is primitive. 
For any element  $g=\left(\begin{smallmatrix} a&b\\c&d \end{smallmatrix}\right)$, let $g'(z)=(cz+d)^{-2}$ and $\chi(g)=\chi(d)$. In this paper, $\chi$ will always be an even character.

\begin{ack}
This work is based on Chapter\,1 of the author's PhD thesis \cite{lou:these}. I would like to sincerely thank my advisors Samuel J. Patterson and Philippe Michel for their encouragements and support. Part of this work has been realized at the University Montpellier 2 and at the Swiss Federal Institute of Technology. I would also like to thank Valentin Blomer for interesting suggestions and advice, and the anonymous referee for a very careful reading of the manuscript.
\end{ack}

\section{Modular forms of half-integral weight}\label{sec:mod}

As usual, $\mathbb{H}=\{z\in\C\,:\, \Im(z)>0\}$ is the upper half-plane and $\Gamma_0(M)$ is the congruence subgroup modulo $M$. Let $\kappa:\Gamma_0(4)\too \{\pm1\}$ be the multiplicative system for the group $\Gamma_0(4)$ defined by

\begin{equation}\label{mod:eq:1}
\gamma'(z)^{1/4} \vartheta(\gamma(z)) = \kappa(\gamma) \vartheta(z), \forall \gamma\in\Gamma_0(4), \forall z\in\mathbb{H},
\end{equation}
where $\vartheta(z) = \sum_{n\in\Z} e(n^2z)$. The Jacobi symbol $(u/v)$ is defined, for $v$ odd, as extension of the Legendre symbol; in particular, $(u/v)=(u/-v)$. The symbol $\kappa$ satisfies $\kappa(\gamma)=\kappa(-\gamma)$, and, if $\gamma=\left(\begin{smallmatrix} a&b\\c&d\end{smallmatrix}\right)\in\Gamma_0(4)$ with $d>0$, then, with the convention that $\ec{arg}(z)\in[0,2\pi[$ for $z\in \C$, we have

\begin{equation}\label{theta:eq:kappa}
\kappa(\gamma)=\begin{cases}\displaystyle
\left(\frac{b}{d}\right)_{\!\!2} \, \varepsilon_d
\begin{cases} i & \ec{if }c>0 \\ 1& \ec{if } c \le 0 \end{cases}  
& \ec{for $c$ even and $b \neq 0$}\\
\begin{cases} i & \ec{if }c>0 \\ 1& \ec{if } c \le 0 \end{cases}  
& \ec{for $c$ even and $b = 0$}.\end{cases}
\end{equation}
Formula \eqref{theta:eq:kappa} has been proved in several places, but with different choices of notations (see e.g. \cite{kub} or \cite[(1.9)-(1.10)]{shi}).
Let $\alpha$ be the multiplier factor of weight $1/2$, i.e.

\begin{equation*}
(gh)'(z)^{1/4} = g'\left(h(z)\right)^{1/4} h'(z)^{1/4} \alpha(g,h), \qquad \forall g,h \in SL_2(\R).
\end{equation*}
Let $N$ be a positive integer. Let $\chi$ be an even Dirichlet character modulo $N$, of conductor $f$ dividing $N$. In particular, the symbol $\kappa$ satisfies

\begin{equation}\label{mod:eq:2}
\kappa\chi(gh) = \kappa\chi(g) \kappa\chi(h) \alpha(g,h), \qquad \forall g,h \in \Gamma_0(4N).
\end{equation}
The cusps of $\Gamma_0(4N)$ are of the form $\sigma^{-1}(\infty)$ with $\sigma\in SL_2(\Z)$. There exists some positive integer $q_\sigma$ such that the subgroup $\Gamma_\sigma=\{\gamma\in \Gamma_0(4N)\,:\, \gamma(\sigma^{-1}(\infty))= \sigma^{-1}(\infty)\}$ is of the form

\begin{equation*}
\Gamma_\sigma= \sigma^{-1} \begin{pmatrix} \pm1& \Lambda_\sigma\\ 0&\pm1 \end{pmatrix} \sigma ,
\end{equation*}
with $\Lambda_\sigma =q_\sigma \Z$. We define, for each cusp $\sigma^{-1}(\infty)$ of $\Gamma_0(4N)$, the real number $\varkappa_\sigma\in[0,1[$, as 

\begin{equation*}
\kappa\chi \left( \sigma\begin{pmatrix}1&q_\sigma \\0&1\end{pmatrix}\sigma\right)=e(-\varkappa_\sigma).
\end{equation*}

\begin{lemme}\label{mod:lemme:1}
Let $\sigma^{-1}(\infty)$ be a cusp of $\Gamma_0(4N)$, and let $\gamma_\sigma=\sigma^{-1} \left(\begin{smallmatrix} \varepsilon & q_\sigma n \\ 0&\varepsilon \end{smallmatrix}\right)\sigma\in\Gamma_\sigma$, for some $\varepsilon=\pm 1$ and $n\in \Z$. Then,

\begin{equation*}
\kappa\chi(\gamma_\sigma) \alpha(\sigma,\gamma_\sigma) = \kappa\chi(\gamma_\sigma) \alpha(\gamma_\sigma,\sigma^{-1})= e(-\varepsilon n \varkappa_\sigma).
\end{equation*}
\end{lemme}
Let $M(N,\chi)$ be the space of functions $f$ holomorphic on $\mathbb{H}$ and at the cusps of $\Gamma_0(4N)$, which satisfy

\begin{equation*}
\gamma'(z)^{1/4} f\left(\gamma(z)\right) =\kappa \chi (\gamma) f(z) \qquad \forall z\in \mathbb{H}, \gamma \in \Gamma_0(4N).
\end{equation*}
With these notations, modular forms have a Fourier expansion of the form

\begin{equation*}
\sigma^{-1}(z)^{1/4} f\left(\sigma(z)\right) = \sum_{n\in \Z} a_f(\sigma,n) e\left(\frac{z(n-\varkappa_\sigma)}{q_\sigma}\right),
\end{equation*}
with $a_f(\sigma,n)=0$ if $n-\varkappa_\sigma<0$. 
The Petersson scalar product on the space $M(N,\chi)$ of modular forms of weight $1/2$ is given by

\begin{equation*}
\langle f,g\rangle = \int_{\Gamma_0(4N)\backslash \mathbb{H}} f(z) \barre{g(z)} \Im(z)^{1/2}\, d\mu(z),
\end{equation*}
where $d\mu (z) =y^{-2}\, dx dy$ is the invariant measure. 

For a character $\psi$, one defines

\begin{equation*}
\vartheta_\psi(z)=\sum \psi(n)e(n^2z).
\end{equation*}
We introduce an other twist of $\vartheta(z)$, namely

\begin{equation}\label{mod:eq:3}
\vartheta_{d,s,q}(z) = \sum_{n\in \Z} \chi\chi_d(n) c_q(n) e(dn^2s^2 z).
\end{equation}
Here, $c_q(n)$ is the Ramanujan's sum, of which one of the representations is 

\begin{equation}\label{mod:eq:ram}
c_q(n)=\sum_{d\tq (q,n)} \mu\left(\frac{q}{d}\right) d.
\end{equation}
The functions $\vartheta_{d,s,q}(z)$ can be expressed in terms of $\vartheta_\psi(z)$ by means of the following equality, valid for any character $\psi$:

\begin{equation}\label{mod:eq:4}
\sum_{n\in\Z}\psi(n) c_q(n) e(n^2z) = \sum_{j\tq q} \mu\left(\frac{q}{j}\right) j \psi(j) \vartheta_\psi(j^2 z).
\end{equation}

Recall that $f_d$ is defined as the conductor of $\chi\chi_d$.
For any triple $(d,s,q)$ of positive integers, consider the condition $(C)$ given by

\begin{equation*}
(C) \qquad \left\{
\begin{aligned}
&d f_d^2 s^2 q^2 \mid N\\
&\ec{$d$ square-free, $s$ supported by $f_d$ and $q$ coprime to $f_d$}.
\end{aligned}\right.
\end{equation*}  

\begin{thm}\label{mod:thm}
Let $N\in \N$. Let $\chi$ be an even character modulo $N$. An orthogonal basis for $M(N,\chi)$ is given by the set 

\begin{equation*}
B= \{\vartheta_{d,s,q} \,:\, (d,s,q) \emph{ satisfies } (C)\}.
\end{equation*}
Moreover,

\begin{equation*}
\norme \vartheta_{d,s,q} \norme^2 
= 2 \pi N \frac{\varphi(q)}{s\sqrt{d}}  \prod_{p\tq f_d}(1-p^{-1}) \prod_{p\tq N}(1+p^{-1}).  
\end{equation*}
\end{thm}

\proof
In order to simplify notations, let us define

\begin{equation*}
\alpha(d,s) = \frac{2 \pi N}{s\sqrt{d}}  \prod_{p\tq f_d}(1-p^{-1}) \prod_{p\tq N}(1+p^{-1}).
\end{equation*}
We know from \cite[Theorem\,A]{ser-sta:theta} that a basis for $M(N,\chi)$ is given by the set 

\begin{equation*}
B_1=\{ \vartheta_\psi(tz)\,:\, \psi \ec{ primitive }; \ec{cond}(\psi)^2 t \tq N; \psi\chi_t (n)=\chi(n), \forall (n,N)=1\}.
\end{equation*} 
By decomposing $t=ds^2q^2$, with $d$ square-free, $s$ supported by $f_d$ and $q$ coprime to $f_d$, one has $\psi=\chi\chi_t=\chi\chi_d$. Therefore, the basis $B_1$ can be written as

\begin{equation*}
B_2=\{ \vartheta_{\chi\chi_d}(ds^2q^2 z)\,:\, f_d^2 ds^2q^2 \mid N\}.
\end{equation*}
We now compute the scalar product of two elements of $B_2$, by using the Rankin-Selberg formula. Let $f$ and $g$ be two elements of $M(N,\chi)$. Denote by $a_f(n)$ and $a_g(n)$ the coefficients of their Fourier expansion at infinity. Let $E(z,s)$ be the Eisenstein series of weight zero for $\Gamma_0(4N)$ defined at the cusp $\infty$. Then, by the Rankin-Selberg method, the following equality holds.

\begin{equation}\label{mod:eq:thm:2}
\langle f,g\rangle = 
\pi N \prod_{p\tq N}(1+p^{-1}) \ec{Res}_{s=1} \left( \sum_{n\geqslant 1} \frac{a_f(n) \barre{a_g(n)}}{n^{s-1/2}}\right).
\end{equation}
Recall that since $N$ is odd, al the primes involved in \eqref{mod:eq:thm:2}, or in the definition of $\alpha(d,s)$, are distinct from 2.
Note that \eqref{mod:eq:thm:2} requires to deal with modular forms, which may be not cuspidal. The Rankin-Selberg method for modular forms of not rapid decay as been developped by Zagier in \cite{zag:rs}, where the author also gives several applications of his main theorem, but none including formula \eqref{mod:eq:thm:2}. To derive formula \eqref{mod:eq:thm:2} from the main theorem of \cite{zag:rs} requires some extra work, however, since a complete proof of \eqref{mod:eq:thm:2} has been given in \cite[Theorem\,2.2]{chi:non-decay}, we can dispense with the details\footnote{The proof of \cite[Theorem\,2.2]{chi:non-decay} contains two errors that luckily neutralize each other. First, in the 5th display of \cite[p.16]{chi:non-decay}, there should be no factor 2 in front of the middle integral. Second, in the 4th display of \cite[p.17]{chi:non-decay}, the right hand side should be multiplied by 2.}.

We apply \eqref{mod:eq:thm:2} with $f=\vartheta_{\chi\chi_d}(ds^2q^2 z)$ 
and $g= \vartheta_{\chi\chi_{d'}}(d's'^2q'^2 z)$. Then, for any $n\geqslant 1$, one has $a_f(n)= 2 \chi\chi_d(m)$, if $n=ds^2q^2m^2$, and $a_g(n)= 2 \chi\chi_{d'}(m')$, if $n=d's'^2q'^2m'^2$. Therefore, non-trivial contributions will occur only if $d=d'$, in which case \eqref{mod:eq:thm:2} gives

\begin{equation}\label{mod:eq:thm:3}
\begin{split}
&\langle \vartheta_{\chi\chi_d}(ds^2q^2 z),\vartheta_{\chi\chi_d}(ds'^2q'^2 z)\rangle\\
&= 4 \pi N \prod_{p\tq N}(1+p^{-1}) 
\times \ec{Res}_{s=1} \Bigg(\sum_{\substack{n\geqslant 1\\n=ds^2q^2m^2\\n=ds'^2q'^2m'^2}} \frac{\chi\chi_d(m) \barre{\chi\chi_{d}(m')}}{n^{s-1/2}}\Bigg).
\end{split}
\end{equation}
Non-trivial contributions will occur for integers $n$ of the form $n=ds^2q^2m^2=ds'^2q'^2m'^2$, with $(m,f_d)=(m',f_d)=1$. Since $s$ and $s'$ are supported by $f_d$ and since $q$ and $q'$ are coprime to $f_d$, this means that $s=s'$ and $n=ds^2q^2q'^2m^2/g^2$, where $g=\ec{gcd}(q,q')$. Then, after further simplifications, we obtain from \eqref{mod:eq:thm:3} that

\begin{equation}\label{mod:eq:thm:4}
\langle \vartheta_{\chi\chi_d}(ds^2q^2 z),\vartheta_{\chi\chi_d}(ds^2q'^2 z)\rangle 
= (q,q') \frac{\barre{\chi\chi_d}(q)}{q} \frac{\chi\chi_d(q')}{q'} \alpha(d,s).
\end{equation}
The orthogonalization of the functions $\vartheta_{\chi\chi_d}(ds^2q^2 z)$, with $q$ varying, is based on the following lemma.

\begin{lemme}\label{mod:lemme:2}
Let $U\in \N$. Let $V$ be the finite dimensional $\C$-vector subspace of a Hilbert space, with scalar product $\langle\pt,\pt\rangle$. Assume that a basis of $V$ is given by the set $E=\{f_u\,:\,u\mid U\}$. Assume that there exists a function $g$ such that $g(d)>0$ and

\begin{equation}\label{mod:eq:lemme2:1}
\langle f_u,f_v\rangle = \sum_{d\tq (u,v)} g(d), \quad \forall u,v \mid U.
\end{equation}
For $u \mid U$, let $f_u' = \sum_{j\tq u} \mu(u/j) f_j$.
Then the set $E'=\{f_u'\,:\,u \mid U\}$ is an orthogonal basis of $V$. Moreover,

\begin{equation}\label{mod:eq:lemme2:2}
\langle f_u',f_v'\rangle =\begin{cases}
g(u)&\ec{if }u=v\\
0&\ec{if }u\neq v.\end{cases}
\end{equation}
\end{lemme}

\proof
It is clear that $E'\subseteq V$ and that \eqref{mod:eq:lemme2:2} follows from the definition of $f_u'$ that

\begin{equation}\label{mod:eq:lemme2:2bis}
\langle f_u',f_v'\rangle =\sum_{j\tq u}\sum_{k\tq v}\mu\left(\frac{u}{j}\right) \mu\left(\frac{v}{k}\right) \langle f_u,f_v\rangle, \quad \forall u,v \mid U.
\end{equation} 
We shall make use of the following M\"obius inversion formula in two variables: 

\begin{equation}\label{mod:eq:lemme2:3}
F(u,v)=\sum_{d\tq u} \sum_{e\tq v} G(d,e) \iff G(u,v)=\sum_{j\tq u}\sum_{k\mid v} \mu\left(\frac{u}{j}\right) \mu\left(\frac{v}{k}\right) F(j,k),
\end{equation}
for any two functions $F,G:\N\times\N\to \R$. Define $F(u,v)=\langle f_u,f_v\rangle$ and $G(u,v)=g(u)$ if $u=v$ and $G(u,v)=0$ if $u\neq v$. Then, by \eqref{mod:eq:lemme2:1}, the first equality of \eqref{mod:eq:lemme2:3} is verified. Thus the second one also holds, which means, by \eqref{mod:eq:lemme2:2bis}, that the vectors $f_u'$ are non-zero orthogonal vectors whose scalar product is given by \eqref{mod:eq:lemme2:2}.\\
\qed

After suitable normalization, and noticing that $(u,v)=\sum_{d\tq (u,v)}\varphi(d)$, one is led to apply Lemma\,\ref{mod:lemme:2} with $f_u(z)= \alpha(d,s)^{-1/2} u \chi\chi_d(u) \vartheta_{\chi\chi_d}(ds^2u^2z)$ and $g(j)=\varphi(j)$. Using formula \eqref{mod:eq:ram}, one sees that the resulting functions $f_u'$ are given by $f_u'=\alpha(d,s)^{-1/2} \vartheta_{d,s,u}(z)$. This concludes the proof of Theorem\,\ref{mod:thm}.\\ 
\qed

Recall that the symbol $\epsilon_c$ is defined, for odd integers $c$, by $\epsilon_c=1$ if $c\congru 1 \Mmod 4)$, and $\epsilon_c=i$ if $c\congru 3 \Mmod 4)$. If $d$ is an odd square-free integer, then $\chi_d$ is a character of conductor $d$ or $4d$, according to if $d\congru 1 \Mmod 4)$ or $d\congru 3 \Mmod 4)$; more precisely,

\begin{equation*}
\chi_d(n)=\left(\frac{n}{d}\right),
\end{equation*}
if $d\congru 1 \Mmod 4)$, and 

\begin{equation*}
\chi_d(n)=\begin{cases}
0&\ec{if }(4d,n)\neq 1\\
\epsilon_n^2 \left(\frac{n}{d}\right)&\ec{if }(4d,n)= 1,
\end{cases}
\end{equation*}
if $d\congru 3 \Mmod 4)$. 

\begin{cor}\label{mod:cor}
Let $D,f$ be odd positive integers coprime to each other. Let $\chi$ be an even primitive character of conductor $f$. 
\begin{itemize}
\item[$(i)$] The space $M(Df,\chi)$ is non-trivial only if $f$ is square-free, $f\congru 1 \Mmod 4)$ and $\chi=\chi_f$.
\item[$(ii)$] If the space $M(Df,\chi)$ is non-trivial, an orthonormal basis of $M(Df,\chi)$ is given by the set $\{\vartheta_{tf,s,q}\}$, where $t^3s^2q^2 \tq D$, $t$ is square-free, $(t,f)=1$, $t\congru 1 \Mmod 4)$, $s$ is supported by $t$, $(q,t)=1$. 
\end{itemize}
\end{cor}

\proof 
Assume that $M(Df,\chi)\neq \{0\}$. By Theorem\,\ref{mod:thm} there exists some element $\vartheta_{d,s,q}\in M(Df,\chi)$, with $(d,s,q)$ satisfying condition $(C)$; in particular, $d$ is an odd square-free integer such that 

\begin{equation}\label{mod:eq:cor}
f_d^2 d\mid Df.
\end{equation}

For $(i)$. Let $p$ be a prime divisor of $f$ and let $\psi$ be the $p$-component of $\chi$, say of order $p^e$, where $\ec{ord}_p(f)=e$. Since $\chi$ is primitive, $e \geqslant 1$. Then $f_d$, the conductor of $\chi\chi_d$ is divisible by the conductor of $\psi\chi_d$. 
Assume that $p\nmid d$ or that $\psi\neq (\pt/p)$. Then the conductor of $\psi\chi_d$ is divisible by $p^e$. Thus $p^e\tq f_d$ and, by \eqref{mod:eq:cor}, one obtains $p^{2e}\tq Df$, which contradicts the fact that $\ec{ord}_p(Df)=\ec{ord}_p(f)=e$. 
Thus $p\tq d$ and $\psi=(\pt/p)$. This shows that $f$ has to be a square-free integer dividing $d$. Moreover, since $\chi$ is even, we conclude that $\chi=\chi_f$, with $f\congru 1 \Mmod 4)$. This proves $(i)$. 

For $(ii)$. Let $t$ be the square-free integer (necessarily coprime to $f$) such that $d=ft$. Then $\chi\chi_d$ is the primitive character associated to $\chi_f\chi_{ft}$. From \eqref{mod:eq:cor}, one sees that $f_d$ has to be odd. This implies $ft\congru 1 \Mmod 4)$, and therefore $t\congru 1 \Mmod 4)$. It follows that $\chi\chi_d=\chi_t$ and that $f_d=t$. The condition $(C)$ then translates into the condition given in $(ii)$.
\qed

\section{The Kloosterman-Selberg Zeta function}\label{sec:zeta}

In this section, we shortly describe the geometric Kloosterman sums, their associated Zeta function, and we obtain Theorem\,\ref{zeta:thm}, as a special case of the Kuznetsov formula. Some more notations have to be introduced: for a matrix $\gamma=\left(\begin{smallmatrix} a&b \\ c&d \end{smallmatrix}\right)$, let $a(\gamma)=a$, $b(\gamma)=b$, $c(\gamma)=c$ and $d(\gamma)=d$. 

\begin{defi}\label{zeta:defi:kloo} 
Let  $\sigma^{-1}(\infty)$ and $\tau^{-1}(\infty)$ be two cusps of $\Gamma_0(4N)$. Let $m,n\in \Z-\{0\}$. Then, for any $c\in\Z$, the geometric Kloosterman sum is defined for positive integers $c$ by

\begin{equation*}
\begin{split}
K_{\sigma,\tau}(m,n;c)&=
\sum_{\substack{\gamma \in \Gamma_\sigma \backslash \Gamma / \Gamma_\tau \\ \tq c(\sigma \gamma\tau^{-1})\tq = c}} 
\barre{\chi \kappa (\gamma) \alpha (\sigma,\gamma) \alpha (\sigma \gamma,\tau^{-1})}\\
& e\left(\frac{(m-\varkappa_\sigma)}{q_\sigma}\frac{a(\sigma \gamma\tau^{-1})}{c(\sigma \gamma\tau^{-1})}\right) e\left(\frac{(n-\varkappa_\tau)}{q_\tau}\frac{d(\sigma \gamma\tau^{-1})}{c(\sigma \gamma\tau^{-1})}\right).
\end{split}
\end{equation*}
\end{defi}
One verifies by  using Lemma\,\ref{mod:lemme:1} the the geometric Kloosterman sums are well defined.
We associate to the sums $K_{\sigma,\tau}(m,n;c)$ the Kloosterman-Selberg Zeta function $Z_{\sigma,\tau,m,n}(s)$ defined by

\begin{equation}
Z_{\sigma,\tau,m,n}(s) =
\sum_{c>0}\frac{K_{\sigma,\tau}(m,n;c)}{c^s}.
\end{equation}

Although an individual bound similar to \eqref{intro:eq:bound} holds for $K_{\sigma,\tau}(m,n;c)$,
we shall not go into details here, but only remark (see the first display in the proof of Lemma\,\ref{zeta:lemme}) that $\tq K_{\sigma,\tau}(m,n;c)\tq \le q_\sigma q_\tau \tq c\tq^2$. This shows that $Z_{\sigma,\tau,m,n}(s)$ is well defined, for $\Re(s)>3$. 

\begin{thm}\label{zeta:thm}
Let $\chi, \sigma,\tau$, $m,n$ be as in Definition\,\ref{zeta:defi:kloo}. Let $X>1$. Then for any $\varepsilon>0$,

\begin{equation*}
\sum_{0<c < X} \frac{K_{\sigma,\tau}(m,n;c)}{c^{1-u}} =  \frac{1}{1+2u} c(N,\chi;\sigma,\tau,m,n) X^{1/2+u} +  \rond{O}\left( X^{1/4+u+\varepsilon}\right),
\end{equation*}
with $c(N,\chi;\sigma,\tau,m,n)=0$ if $m-\varkappa_\sigma<0$ or $n-\varkappa_\sigma<0$, and 
otherwise

\begin{equation*}
c(N,\chi;\sigma,\tau,m,n)= \frac{2(1+i)}{\pi} q_\sigma q_\tau \sum_{f\in B(N,\chi)} \barre{a_f(\sigma,m)} \ a_f(\tau,n),
\end{equation*}
where $B(N,\chi)$ is any orthonormal basis of the space $M(N,\chi)$, defined in Section\,\ref{sec:mod}.
\end{thm}

\proof
The argument is taken from \cite{gol-sar:kloo}: the theory of Poincar\'e series allows us to continue meromorphically $Z_{\sigma,\tau,m,m}(s)$ to $\Re(s) >1$ and shows that its poles are located at $s=2s_i$, where the $s_i(1-s_i)$ are the exceptional eigenvalues of the hyperbolic Laplacian of weight $1/2$.  We write the spectral parameters $s_i$ as $s_1=3/4 > s_2 > \pts > 1/2$. The Laplacian operator is

\begin{equation*}
-y^2 \left(\frac{\partial^2}{\partial x^2}+\frac{\partial^2}{\partial y^2}\right) +\frac{iy}{2} \frac{\partial }{\partial x}.
\end{equation*}
The space $L^2\left(\Gamma_0(4N)\backslash \mathbb{H}, 1/2, \chi\kappa, \lambda_i\right)$ is the space of Maa\ss{} forms, i.e. of eigenfunctions of the Laplacian, having polynomial growth at the cusps, and satisfying 

\begin{equation*}
j_\gamma(z) f\left(\gamma(z)\right) =\kappa \chi (\gamma) f(z) \qquad \forall z\in \mathbb{H}, \gamma \in \Gamma_0(4N).
\end{equation*}
The eigenvalues can be written $\lambda=s(1-s)$ with $\Re(s)\geqslant 1/2$. We denote by $L^2(N,\chi,s)$ the corresponding subspace of $L^2\left(\Gamma_0(4N)\backslash \mathbb{H}, 1/2, \chi\kappa, \lambda_i\right)$.
The Fourier expansion of $f\in L^2(N,\chi,s)$ at a cusp $\sigma^{-1}(\infty)$ is given by 

\begin{equation*}
\begin{split}
j_{\sigma^{-1}}(z) f(\sigma^{-1}(z)) 
&= \delta_{s>1/2}\delta_{\varkappa_\sigma=0} \rho_f(\sigma,0)y^{1-s} \\
&+ \delta_{\varkappa_\sigma\neq 0} \rho_f(\sigma,0) W_{\frac{-1}{4},s-\frac{1}{2}} \left(4\pi \frac{\varkappa_\sigma}{q_\sigma} y\right) e\left( \frac{-\varkappa_\sigma}{q_\sigma}x\right)\\
&+ \sum_{\substack{n\in \Z\\ n\neq 0}}\rho_f(\sigma,n) W_{\frac{\ec{sgn}(n)}{4}, s-1/2}\left(4\pi \frac{\tq n-\varkappa_\sigma\tq}{q_\sigma}y\right) e\left(\frac{n-\varkappa_\sigma}{q_\sigma}x\right).
\end{split}
\end{equation*}
One can derive an asymptotic formula for the $K_{\sigma,\tau}(m,n;c)$ from the analytic properties of its Zeta function $Z_{\sigma,\tau,m,n}(s)$. Let $\chi, \sigma,\tau$, $m,n$ as above. Let $X>1$. Then for any $\varepsilon>0$,

\begin{equation*}
\sum_{0<c < X} \frac{K_{\sigma,\tau}(m,n;c)}{c} =  2 \ec{Res}_{s=3/2}\left( Z_{\sigma,\tau,m,n}(s)\right) X^{1/2} + \rond{O}\left( X^{1/4+\varepsilon}\right).
\end{equation*}
The residues of $Z_{\sigma,\tau,m,n}(s)$ can be expressed as follows:  

\begin{equation*}
\begin{split}
\ec{Res}_{s=2s_i} \left(Z_{\sigma,\tau,m,n}(s)\right)
= e^{i\pi/4} \frac{4^{1-s_i}}{\pi^{2s_i-1/2}}  q_\sigma q_\tau \left(\left\vert n-\frac{\varkappa_\tau}{q_\tau}\right\vert \left\vert m-\frac{\varkappa_\sigma}{q_\sigma}\right\vert \right)^{1-s_i}\\
\times \Gamma(2s_i-1) \frac{\Gamma(s_i + \ec{sgn}(n)/4)}{\Gamma(s_i - \ec{sgn}(m)/4)} \sum_{u\in \rond{B}(s_i)} \barre{\rho_u(\sigma,m)} \rho_u(\tau,n),
\end{split}
\end{equation*}
where $\rond{B}(s_i)$ is an orthonormal basis of $L^2(N,\chi,s_i)=L^2\left(\Gamma_0(4N)\backslash \mathbb{H}, 1/2, \chi\kappa, \lambda_i\right)$.
There exists an isomorphism $M(N,\chi)\rightarrow L^2(N,\chi,3/4)$, given by $f(z)\mapsto u(z)=f(z)\Im(z)^{1/4}$; the injectivity of the homomorphism is clear since modular forms of weight $1/2$ are square-integrable, and the surjectivity comes from \cite[Satz\,9.1]{roe:autoformen}. Under this isomorphism, the Fourier coefficients satisfy, for $n >0$, the relation $\rho_u(\sigma,n)= a_f(\sigma,n) \left(4\pi \frac{n-\varkappa_\sigma}{q_\sigma}\right)^{-1/4}$. This concludes the proof of Theorem\,\ref{zeta:thm}.
\qed

We conclude this section with a lemma, giving the connection between the geometric and the arithmetic Sali\'e sums. 

\begin{lemme}\label{zeta:lemme}
Let $\Gamma=\Gamma_0(4Df)$ with $D$ and $f$ mutually coprime odd integers. Let $\sigma^{-1}(\infty)$ and $\tau^{-1}(\infty)$ be the two cusps of $\Gamma$ defined by $\tau^{-1}=Id$ and

\begin{equation*}
\sigma^{-1} =  \begin{pmatrix} \alpha&\beta\\ D&4 f \end{pmatrix}  \in SL_2(\Z).
\end{equation*}
Then,

\begin{gather*}
K_{\sigma,\tau}(4fm,n;c)=\\
\begin{cases}
0& \ec{ if $c\notcongru 0\Mmod D)$ or if $c$ is even}\\
\barre{\epsilon_D} \barre{\epsilon_c} \chi(c) \barre{\chi}(D)  \left(\frac{f}{D}\right) K_2(m,n;c) & \ec{ if $c\congru 0\Mmod D)$ and if $c$ is odd}.
\end{cases}
\end{gather*}
\end{lemme}
 
\proof
With our assumptions, the definition of $K_{\sigma,\tau}(m,n;c)$ can be written as

\begin{equation*}
\begin{split}
K_{\sigma,\tau}(m,n;c) 
&= \sum_{\substack{ a\mmod q_\sigma c)\\d \mmod q_\tau c)\\ \sigma^{-1} \left(\begin{smallmatrix} a&*\\c&d\end{smallmatrix}\right)\tau \in \Gamma}} 
\barre{\kappa\chi}\left(\sigma^{-1}\begin{pmatrix} a&*\\c&d\end{pmatrix}\tau\right) \barre{\alpha}\left(\sigma,\sigma^{-1}\begin{pmatrix} a&*\\c&d\end{pmatrix}\tau\right) \\ 
&\times e\left(\frac{(m-\varkappa_\sigma)a}{cq_\sigma}\right) e\left(\frac{(n-\varkappa_\tau) d}{cq_\tau}\right).
\end{split}
\end{equation*}

The condition appearing in the sum means that $ad \congru 1\Mmod c)$ et $c(\sigma^{-1} \left(\begin{smallmatrix} a&*\\c&d\end{smallmatrix}\right)\tau) \congru 0\Mmod 4fD)$. With our choice of $\sigma$ and $\tau$, we obtain the conditions $c \congru 0 \Mmod D)$ and $a\congru 0 \Mmod 4f)$. Then,

\begin{align*}
&\barre{\kappa \chi} \left( \sigma^{-1} \begin{pmatrix} a&b \\ c&d \end{pmatrix}\right) \barre{\alpha}\left(\sigma, \sigma^{-1} \begin{pmatrix} a&b \\ c&d \end{pmatrix}\right)\\
&= \barre{\kappa \chi} \left( \sigma^{-1}S (-S) \begin{pmatrix} a&b \\ c&d \end{pmatrix}\right) \barre{\alpha}\left(\sigma, \sigma^{-1}S (-S) \begin{pmatrix} a&b \\ c&d \end{pmatrix}\right)\\
&= \barre{\kappa \chi}\left( \sigma^{-1}S\right) \barre{\kappa \chi}\left(-S \begin{pmatrix} a&b \\ c&d \end{pmatrix}\right) 
\barre{\alpha}\left(\sigma^{-1}S,-S \begin{pmatrix} a&b \\ c&d \end{pmatrix} \right)
\barre{\alpha}\left(\sigma, \sigma^{-1}S (-S) \begin{pmatrix} a&b \\ c&d \end{pmatrix}\right)\\
&= \barre{\kappa \chi}\left( \sigma^{-1}S\right) \barre{\kappa \chi}\left(-S \begin{pmatrix} a&b \\ c&d \end{pmatrix}\right)\
\barre{\alpha}\left(\sigma , \sigma^{-1}S\right) \barre{\alpha}\left(S, -S \begin{pmatrix} a&b \\ c&d \end{pmatrix} \right).
\end{align*}
Let us assume, as we may, that $a>0$ et $d>0$. Then $\barre{\alpha}\left(\sigma , \sigma^{-1}S\right) \barre{\alpha}\left(S, -S \left(\begin{smallmatrix} a&b \\ c&d \end{smallmatrix}\right) \right)=1$. 
Since $b>0$, we obtain from the definition of $\kappa$ that
 
\begin{align*}
\kappa\left(-S \begin{pmatrix} a&b \\ c&d \end{pmatrix}\right) 
&=\kappa\left(\begin{pmatrix} c&d \\ -a&-b \end{pmatrix}\right)=\kappa\left(\begin{pmatrix} -c&-d \\ a&b \end{pmatrix}\right)\\
&= \left(\frac{-d}{b}\right) \epsilon_b i=\left(\frac{-1}{b}\right)\left(\frac{d}{b}\right) \epsilon_b i .
\end{align*}
Since $bc\congru -1\Mmod 4)$, we have $\epsilon_{b} (-1/b) i =\epsilon_c$. It remains

\begin{equation*}
\kappa\left(-S \begin{pmatrix} a&b \\ c&d \end{pmatrix}\right) 
= \left(\frac{d}{b}\right) \epsilon_c= \left(\frac{a}{b}\right) \epsilon_c = \left(\frac{a}{c}\right) \epsilon_c.
\end{equation*}
One computes also

\begin{equation*}
\kappa(\sigma^{-1}S) = \kappa \left( \begin{pmatrix} \beta & -\alpha \\ 4f & -D \end{pmatrix}\right) =\kappa \left( \begin{pmatrix} -\beta & \alpha \\ -4f & D \end{pmatrix}\right)= \left(\frac{\alpha}{D}\right) \epsilon_D,
\end{equation*}
as well as

\begin{equation*}
\chi(\sigma^{-1}S) =\chi(-D) \ , \quad \ec{ and }  \chi\left( -S \begin{pmatrix} a&b \\ c&d \end{pmatrix}\right) =\chi(-b).
\end{equation*}
Finally, we showed that, if $c\congru 0\Mmod D)$,

\begin{gather*}
K_{\sigma,\tau}(m,n;c)=
\sum_{\substack{a\mmod q_\sigma c)\\d\mmod q_\tau c)\\ ad\congru 1\mmod c)\\ a\congru 0\mmod 4f)}} 
\barre{\left(\frac{\alpha}{D}\right) \epsilon_D  \left(\frac{a}{c}\right) \epsilon_c \chi(-D) \chi(-b)} e\left(\frac{m a}{q_\sigma c}\right) e\left( \frac{n d}{q_\tau c}\right).
\end{gather*}
In our case, $q_\sigma=4f$ and $q_\tau=1$. This leads, after some further simplifications, to the assertion of Lemma\,\ref{zeta:lemme}.
\qed

\section{Asymptotic behavior of Sali\'e sums}\label{sec:asymp}

In this section we prove Theorem\,\ref{intro:thm}. For it, we use the formula of Theorem\,\ref{zeta:thm}, with $N=Df$ and $\sigma,\tau$ chosen as in Lemma\,\ref{zeta:lemme}. The left hand side of the formula of Theorem\,\ref{zeta:thm} is then determined by Lemma\,\ref{zeta:lemme}. For the right hand side, we use our determination of an orthogonal basis of $M(Df,\chi)$ given in Corollary\,\ref{mod:cor}, i.e. we use the orthonormal basis formed by the $\vartheta_{tf,s,u}(z)$, where $t^3\mid D$, $s^2\mid D/t^3$, $s$ is supported by $t$, $u^2\mid D/t^3$, $u$ is coprime to $t$.  Recall that one such element $\vartheta_{tf,s,u}(z)$ is defined in \eqref{mod:eq:3}. Therefore, it remains to obtain the Fourier expansion of $\vartheta_{tf,s,u}(z)$ at the cusp $\sigma^{-1}(\infty)$ of $\Gamma_0(4Df)$; this is done in the following proposition. 

\begin{prop}\label{asymp:prop}
Let $\sigma$ be as in Lemma\,\ref{zeta:lemme}. Let $f$ be an odd square free integer, $f>0$,
$f\congru 1\Mmod 4)$ and let $\chi=\chi_f$. Let $D$ be odd and coprime to $f$. Then, 

\begin{equation*}
\sigma'(z)^{-\frac{1}{4}} \vartheta_{tf,s,u} (\sigma^{-1}(z))= \sum_{m\geqslant 0} a_{tf,s,u}(\sigma,m) e\left(\frac{mz}{4f}\right),
\end{equation*}
where, for $m>0$, 

\begin{equation*}
a_{tf,s,u}(\sigma,m) =\begin{cases}
0 &\ec{ if $m \not\in ts^2 \Z^2$},\\
\frac{(1+i)}{\sqrt{f}} \epsilon_D \left(\frac{2m'}{t}\right) 
c_u(m') &\ec{ if } m= m'^2 ts^2.
\end{cases}
\end{equation*}
\end{prop}

\proof
The cusp $\sigma^{-1}(\infty)$ is of width $4f$, and $\varkappa_\sigma=0$. From  \eqref{mod:eq:4}, we have 

\begin{equation*}
\vartheta_{tf,s,u} (z) =  \sum_{j \tq u} \mu \left(\frac{u}{j}\right) j \left(\frac{j}{t}\right) \vartheta_{\chi_t} (tf s^2 j^2 z),
\end{equation*}
so that the first step is to study $\vartheta_{\chi_t} (tf s^2 j^2 z)$ at $\sigma^{-1}(\infty)$.
More generally, we have the 
 
\begin{lemme}\label{asymp:lemme}
Let $\psi$ be an even character of conductor $t$. Let $\sigma$ be as in Lemma\,\ref{zeta:lemme}, and let $T$ be such that $t^2T \tq D $. Then, the Fourier expansion of $\theta_\psi(Tz)$ at $\sigma^{-1}(\infty)$ is

\begin{equation*}
(\sigma^{-1})'(z)^{1/4} \vartheta_\psi \left(Tf \sigma^{-1}(z)\right)
= \frac{(1+i)}{2\sqrt{f}} \left(\frac{f\alpha}{DT}\right) \epsilon_{DT}\psi(2) \vartheta_\psi \left(\frac{Tz}{4f}\right).
\end{equation*}
\end{lemme}

We postpone the proof of Lemma\,\ref{asymp:lemme} at the end of this section. As a consequence, with $\psi=\chi_t$ and $T=ts^2j^2$, we have 

\begin{align*}
& (\sigma^{-1})'(z)^{\frac{1}{4}} \vartheta_{tf,s,u} (\sigma^{-1}(z))\\
&= \frac{(1+i)}{2\sqrt{f}} \left(\frac{f\alpha}{Dt}\right) \epsilon_{Dt}\left(\frac{2}{t}\right)\sum_{j \tq u} \mu \left(\frac{u}{j}\right) j \chi_t(j) \vartheta_{\chi_t} \left(\frac{ts^2j^2z}{4f}\right).
\end{align*}
From $4\alpha f \congru 1 \Mmod D)$ follows that $(f\alpha/Dt)=1$. Then, 

\begin{align*}
&\sum_{j \tq u} \mu \left(\frac{u}{j}\right) j \left(\frac{j}{t}\right) 
 \vartheta_{\chi_t} \left(\frac{ts^2j^2z}{4f}\right)\\
&=\sum_{j \tq u} \mu \left(\frac{u}{j}\right) j \left(\frac{j}{t}\right) 
\sum_{m\in \Z} \left(\frac{m}{t}\right) e\left(\frac{ts^2j^2m^2z}{4f}\right)\\
&= \sum_{m'>0} \sum_{j \tq (u,m')} \mu \left(\frac{u}{j}\right) j \left(\frac{j}{t}\right) 
\left( \left(\frac{m'/j}{t}\right)+\left(\frac{-m'/j}{t}\right)\right) e\left(\frac{ts^2m'^2z}{4f}\right)\\
&=  \sum_{m'>0} \left(1+\left(\frac{-1}{t}\right)\right) \left(\frac{m'}{t}\right)c_u(m') e\left(\frac{ts^2m'^2z}{4f}\right),
\end{align*}
since $(u,t)=1$. We conclude by using the fact that $t\congru 1 \Mmod 4)$. 
\qed

\proof[Proof of Theorem\,\ref{intro:thm}]
Let $m,n\in\Z$, $m,n$ positive. Let $\chi=\chi_f$ with $f$ square-free and $f\congru 1 \Mmod 4)$. 
Let 

\begin{equation*}
c(D,f;t,s,q) 
= \left(\frac{2 \pi D \sqrt{f}}{s\sqrt{t}}\varphi(q) \prod_{p\tq t}(1-p^{-1}) \prod_{p\tq Df}(1+p^{-1})\right)^{-1/2}.  
\end{equation*}
By Corollary\,\ref{mod:cor}, an orthonormal basis $B(Df,\chi)$ of $M(Df,\chi)$ is given by 

\begin{equation*}
B(Df,\chi) = \{\vartheta_{t,s,q}'\}
\end{equation*}
where $\vartheta_{t,s,q}'(z)= c(D,f;t,s,q) \vartheta_{tf,s,q}(z)$, and where the parameters satisfy $t\congru 1 \Mmod 4)$, $t$ is square-free and coprime to $f$, $s\mid t^\infty$, $(q,t)=1$ and $t^3s^2q^2 \mid D$.
Let $\tau$ and $\sigma$ be as in Lemma\,\ref{zeta:lemme}. The $n$-th Fourier coefficient of $\vartheta_{t,s,q}'(z)$ at $\tau^{-1}(\infty)$ is 

\begin{equation}\label{asymp:eq:thm:1}
a_{t,s,q}(\tau,n)=\begin{cases} 0&\ec{if $n\not\in tfs^2 \Z^2$},\\
2 c(D,f;t,s,q) \left(\frac{n'}{t}\right) c_q(n')&\ec{if $n=tfs^2n'^2$}.
\end{cases}
\end{equation}
By Proposition\,\ref{asymp:prop}, the $m$-th Fourier coefficient of $\vartheta_{t,s,q}'(z)$ at $\sigma^{-1}(\infty)$ is 

\begin{equation}\label{asymp:eq:thm:2}
a_{t,s,q}(\sigma,m) =\begin{cases}
0 &\ec{ if $m \not\in ts^2 \Z^2$},\\
\frac{(1+i)}{\sqrt{f}} \epsilon_D \left(\frac{2m'}{t}\right) c(D,f;t,s,q)
c_q(m') &\ec{ if } m= ts^2m'^2.
\end{cases}
\end{equation}

Let us fix some $\vartheta_{t,s,q}'\in B(Df,\chi)$ and look at the expression $\barre{a_{tf,s,q}}(\sigma,4fm) a_{tf,s,q}(\tau,n)$. From \eqref{asymp:eq:thm:1} and  \eqref{asymp:eq:thm:2}, we see that this expression is non-zero only if $n\in tfs^2\Z^2$ and $4fm\in t s^2 \Z^2$. Since $(2f,t)=1$ and since $f$ is square-free, the second condition means $m\in tfs^2\Z^2$; this shows that condition $(i)$ of Theorem\,\ref{intro:thm} is necessary. Moreover, writing $m=tfs^2m'^2$ and $n=tfs^2n'^2$, since the factor $(m'n'/t)$ appears in $\barre{a_{tf,s,q}}(\sigma,4fm) a_{tf,s,q}(\tau,n)$, we also obtain condition $(ii)$.  

Thus the couple $(t,s)$ is completely determined by $m$ and $n$, and, writing $m=tfs^2m'^2$ and $n=tfs^2n'^2$, we get 

\begin{equation*}
\barre{a_{tf,s,q}}(\sigma,4fm) a_{tf,s,q}(\tau,n)
= \frac{2(1-i)}{\sqrt{f}} \barre{\epsilon_D} \left(\frac{fm'n'}{t}\right) c(D,f;t,s,q)^2 c_q(n') c_q(m').
\end{equation*}
Therefore, if we decompose $D=D_tD'$, with $D_t$ supported by $t$, and $D'$ coprime to $t$,

\begin{align}
&\sum_{\vartheta_{t,s,q}' \in B(Df,\chi)} \barre{a_{t,s,q}}(\sigma,4fm) a_{t,s,q}(\tau,n)\nonumber\\
&= \frac{2(1-i)}{\sqrt{f}} \barre{\epsilon_D} \left(\frac{fm_1n_1}{t}\right)
\sum_{q^2 \tq D'}c(D,f;t,s,q)^2 c_q(n') c_q(m')\nonumber \\
&= \frac{(1-i)}{\pi} \frac{s\sqrt{t}}{Df} \barre{\epsilon_D} \left(\frac{fm_1n_1}{t}\right) \prod_{p\tq t}(1-p^{-1})^{-1} \prod_{p\tq Df}(1+p^{-1})^{-1} \nonumber\\
& \times \sum_{q^2 \tq D'} \varphi(q)^{-1} c_q(n') c_q(m').\label{asymp:eq:thm:3}
\end{align}
Let us define 

\begin{equation}
T(m,n;c)= \sum_{u\tq c}\varphi(u)^{-1} c_u(n) c_u(m).
\end{equation}

\begin{lemme}\label{ram:lemme:2}
If $(m,c)\neq (n,c)$, then $T(m,n;c)=0$. If $(m,c)= (n,c)$, then 

\begin{gather*}
T(m,n;c) = T(n,m;c) =
(m,c)  \prod_{p\tq \frac{c}{(m,c)}} \frac{p}{p-1}. 
\end{gather*}
\end{lemme}

\proof
Let ud first note that $T(m,n;c)$ is, as a function of $c$, multiplicative. Let us fix some $p\tq c$. For $a\in \N$, let us denote by $e_p(a)$ the order of the prime $p$ in $a$. We may assume that $e_p(m)\leqslant e_p(n)$. Then, from the determination of Ramanujan's sums at prime power, we obtain

\begin{align*}
&T(m,n;p^{e_p(c)}) = \sum_{i=0}^{\min(e_p(c),e_p(m)+1)} \varphi(p^i)^{-1} c_{p^i}(m) c_{p^i}(n)\\
&= \sum_{i=0}^{\min(e_p(c),e_p(m))} \varphi(p^i)+\delta_{e_p(c)>e_p(m)}\frac{p^{e_p(m)}}{p-1} \begin{cases}
1& \ec{if }e_p(n)=e_p(m)\\
1-p& \ec{if }e_p(n)>e_p(m)
\end{cases}\\
&= p^{\min(e_p(c),e_p(m))} +\delta_{e_p(c)>e_p(m)}\frac{p^{e_p(m)}}{p-1} \begin{cases}
1& \ec{if }e_p(n)=e_p(m)\\
1-p& \ec{if }e_p(n)>e_p(m)
\end{cases}\\
&= \begin{cases} 
p^{e_p(c)}&\ec{if } e_p(c)\leqslant e_p(m)\leqslant e_p(n)\\
p^{e_p(m)} \frac{p}{p-1}
&\ec{if } e_p(n)=e_p(m)<e_p(c)\\
0&\ec{if }e_p(m)<e_p(c) \ec{ and } e_p(m)\neq e_p(n).\end{cases}
\end{align*}
Thus for every $p\tq c$, a non-zero contribution occurs only if, either $e_p(m)\geqslant e_p(c)$ and $e_p(n)\geqslant e_p(c)$, or $e_p(m)=e_p(n)< e_p(c)$; these two conditions can be rewritten as $(m,c)=(n,c)$.
\qed

We conclude the proof of Theorem\,\ref{intro:thm} by combining \eqref{asymp:eq:thm:3}, Lemma\,\ref{ram:lemme:2}, Theorem\,\ref{zeta:thm} and Lemma\,\ref{zeta:lemme}. 
\qed

\proof[Proof of Lemma\,\ref{asymp:lemme}]
We shall use a more general theta function. Let $\eta, \eta' \in \R$. Define

\begin{equation}\label{asymp:eq:lemme:1}
\theta\left(\eta, \eta', z\right) = \sum_{n \in \Z}e\left(\frac{1}{2}\left(n+\frac{\eta}{2}\right)^2 z\right) e\left(\left(n+\frac{\eta}{2}\right)\frac{\eta'}{2}\right).
\end{equation}
Then if $\gamma = \left(\begin{smallmatrix}a&b \\ c&d \end{smallmatrix}\right) \in \rm{SL}_2(\Z)$, $z \in \mathbb{H}$ et $\eta,\eta' \in \R$ we have (see \cite{far-kra:theta}, Theorem\,1.11 p. 81):

\begin{equation*}
\theta(\eta,\eta',\gamma(z)) = c(\eta,\eta',\gamma) (cz+d)^{1/2} \theta(a\eta+c\eta'-ac,b\eta+d\eta'+bd,z),
\end{equation*}
with a constant $c(\eta,\eta',\gamma)$ satisfying

\begin{equation*}
c(\eta,\eta',\gamma) = e\left(-\frac{(a\eta+c\eta')bd}{4}-\frac{(ab\eta^2+cd\eta'^2+2bc\eta\eta')}{8}\right) c(0,0,\gamma).
\end{equation*}
The particular case $\eta'=0$ gives

\begin{equation}\label{asymp:eq:lemme:2}
\theta(\eta,0,\gamma(z)) = c(\eta,0,\gamma) (cz+d)^{1/2} \theta(a\eta-ac,b\eta+bd,z),
\end{equation}
with 

\begin{equation}\label{asymp:eq:lemme:3}
c(\eta,0,\gamma) = e\left(-\frac{\eta abd}{4}-\frac{ab\eta^2}{8}\right) c(0,0,\gamma).
\end{equation}
Assume that $\gamma\in\Lambda$; Then $ac$ and $bd$ are even, and we have

\begin{align*}
&\theta(a\eta-ac,b\eta+bd,z)\\
&= \sum_{n \in \Z}e\left(\frac{1}{2}\left(n+\frac{a\eta-ac}{2}\right)^2 z\right) e\left(\frac{1}{2}\left(n+\frac{a\eta-ac}{2}\right)(b\eta+bd)\right)\\
&= \sum_{n \in \Z}e\left(\frac{1}{2}\left(n+\frac{a\eta}{2}\right)^2 z\right) e\left(\frac{1}{2}\left(n+\frac{a\eta}{2}\right)(b\eta+bd)\right)\\
&= e\left(\frac{a\eta}{4}(b\eta+bd)\right) \sum_{n \in \Z}e\left(\frac{1}{2}\left(n+\frac{a\eta}{2}\right)^2 z\right) e\left(\frac{1}{2}nb\eta\right),
\end{align*}
and introducing this in \eqref{asymp:eq:lemme:2} gives, with \eqref{asymp:eq:lemme:3},

\begin{equation}\label{asymp:eq:lemme:4}
\begin{split}
\theta(\eta,0,\gamma(z)) &= e\left(\frac{ab\eta^2}{8}\right) c(0,0,\gamma) (cz+d)^{1/2} \\
&\sum_{n \in \Z}e\left(\frac{1}{2}\left(n+\frac{a\eta}{2}\right)^2 z\right) e\left(\frac{1}{2}nb\eta\right).
\end{split}
\end{equation}
Remark that $\theta(0,0,z)$ is the function $\theta(z)=\sum_{n\in \Z}e^{i\pi n^2 z}$, which is modular for the group $\Lambda$ generated by $\left(\begin{smallmatrix} 1& 2 \\ 0&1 \end{smallmatrix}\right)$ and $\left(\begin{smallmatrix} 0& -1 \\ 1&0 \end{smallmatrix}\right)$. In particular, there exists a function $\kappa_\theta$ on $\Lambda$ such that 

\begin{equation*}
\gamma'(z)^{1/4} \theta\left(\gamma(z)\right) =\kappa_\theta(\gamma) \theta(z) \qquad \forall z\in \mathbb{H},  \forall \gamma \in \Lambda.
\end{equation*}
Because of $\kappa_\theta(-Id)=1$, $\kappa_\theta$ is determined by its values on the elements $\gamma=\left(\begin{smallmatrix}a&b \\ c&d\end{smallmatrix}\right)\in \Lambda$, with $d>0$; for such an element $\gamma$, and with our choice of the branch of $z^{1/2}$, we have

\begin{equation*}
\kappa_\theta(\gamma) =  \begin{cases}
e^\frac{i\pi}{4} \left(\frac{2a}{c}\right) \epsilon_c   
& \ec{for $c$ odd and $a \neq 0$}\\
e^\frac{i\pi}{4}  
& \ec{for $c$ odd and $a = 0$}\\
\left(\frac{2b}{d}\right) \epsilon_d 
\begin{cases} i & \ec{if }c>0 \\ 1& \ec{if } c \le 0 \end{cases}  
& \ec{for $c$ even and $b \neq 0$}\\
\begin{cases} i & \ec{if }c>0 \\ 1& \ec{if } c \le 0 \end{cases}  
& \ec{for $c$ even and $b = 0$}.
\end{cases}
\end{equation*}
Thus the constant $c(0,0,\gamma)$ is defined, in the case $\gamma\in\Lambda$, as

\begin{equation*}
\theta(0,0,\gamma(z)) = c(0,0,\gamma) (cz+d)^{1/2} \theta(0,0,z).
\end{equation*}
Since $\theta(0,0,z)=\theta(z)$, we have

\begin{equation*}
c(0,0,\gamma) (cz+d)^{1/2} = \kappa_\theta (\gamma) \gamma'(z)^{-1/4}.
\end{equation*}
Therefore, we conclude that

\begin{equation}\label{asymp:eq:lemme:5}
\theta(\eta,0,\gamma(z)) = e\left(\frac{ab\eta^2}{8}\right) \kappa_\theta (\gamma) \gamma'(z)^{-1/4} \sum_{n \in \Z}e\left(\frac{1}{2}\left(n+\frac{a\eta}{2}\right)^2 z\right) e\left(\frac{1}{2}nb\eta\right).
\end{equation}

Before starting with the computation of the Fourier expansion of $\vartheta_\psi(Tf z)$ at $\sigma^{-1}(\infty)$, we merely remark that

\begin{equation*}
\tilde{\sigma} = \begin{pmatrix}  2f\alpha &  t^2 T\beta \\ D/t^2 T & 2\end{pmatrix}
\end{equation*}
is an element of $\Lambda$ and satisfies

\begin{equation*}
2 t^2 Tf \sigma^{-1}(z) = \tilde{\sigma} (z'), \qquad \ec{with }  z'=\frac{t^2 T }{2 f} z .
\end{equation*}
Now, 

\begin{align*}
&\vartheta_\psi \left(Tf \sigma^{-1}(z)\right)\\
&= \sum_{n\in\Z} \psi(n) e\left(n^2 Tf \sigma^{-1}(z)\right) 
= \sum_{h  (t)} \psi (h) \sum_{n \in \Z} e\left( \left(n + \frac{h}{t}\right)^2 t^2 Tf \sigma^{-1}(z)\right)\\
&= \sum_{h (t)} \psi (h) \sum_{n \in \Z} e\left( \frac{1}{2}\left(n + \frac{h}{t}\right)^2  \tilde{\sigma}(z')\right)
= \sum_{h (t)} \psi (h) \theta\left(\frac{2h}{t}, 0, \tilde{\sigma}(z')\right)
\end{align*}
and by the formula \eqref{asymp:eq:lemme:5}, since $(2\alpha f,t)=1$, 

\begin{align*}
&\vartheta_\psi \left(Tf \sigma^{-1}(z)\right)\\
&= (\tilde{\sigma})' (z')^{-1/4} \kappa_\theta(\tilde{\sigma}) \sum_{h (t)} \psi (h) \sum_{n \in \Z}  e\left(\left(n+ \frac{2  \alpha h f }{t}\right)^2 \frac{z'}{2} \right) \\
&= (\tilde{\sigma})' (z')^{-1/4} \kappa_\theta(\tilde{\sigma}) \barre{\psi}(2\alpha f) \sum_{h (t)} \psi (2\alpha fh) \sum_{n \in \Z} e\left(\left(nt+ 2 \alpha h f\right)^2 \frac{z'}{2t^2} \right) \\
&= (\tilde{\sigma})' (z')^{-1/4} \kappa_\theta(\tilde{\sigma}) \psi(2) \vartheta_\psi \left(\frac{z'}{2t^2}\right),
\end{align*}
since $4\alpha f \congru 1 \Mmod t)$.
We finally verify that

\begin{equation*}
(\tilde{\sigma})'(z')^{-1/4}=\frac{(\sigma^{-1})'(z)^{-1/4}}{\sqrt{2 f}}  
\end{equation*}
and that

\begin{equation*}
\kappa_\theta(\tilde{\sigma}) =e^{i\pi /4} \left(\frac{f\alpha}{DT}\right) \epsilon_{DT}.
\end{equation*}
This finishes the proof of Lemma\,\ref{asymp:lemme}.
\qed

\section{Numerical examples}\label{sec:num}

\begin{figure}[!ht]
\hspace{-0.5cm}\includegraphics[width=0.5\textwidth]{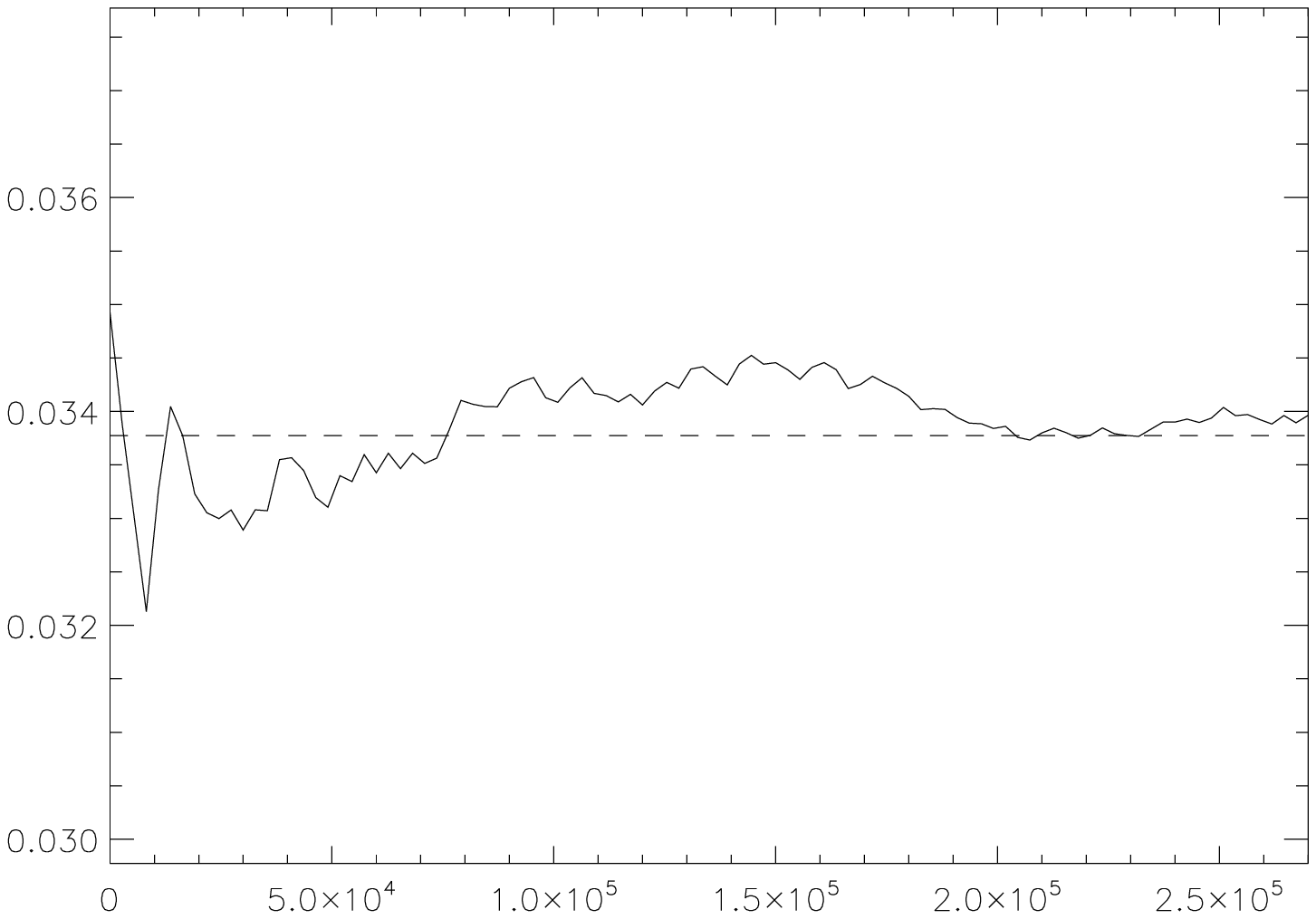}
\hspace{-0.5cm}\includegraphics[width=0.5\textwidth]{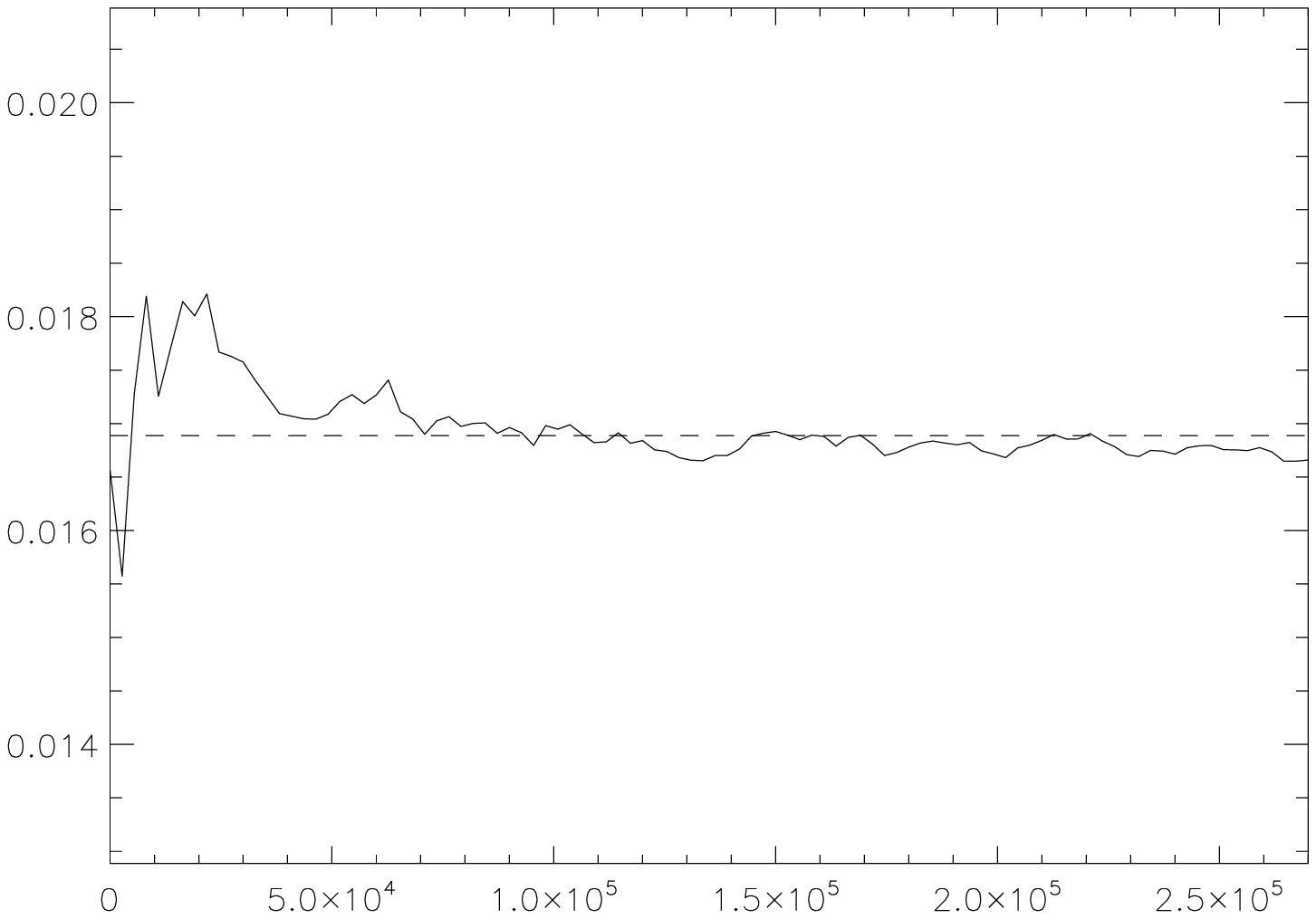}
\caption{{\it Left:} m=n=1, D=27, f=1. \ {\it Right:}  m=n=1, D=45, f=1.}
\label{fig:1}      
\end{figure}

\begin{figure}[!ht]
\hspace{-0.5cm}\includegraphics[width=0.5\textwidth]{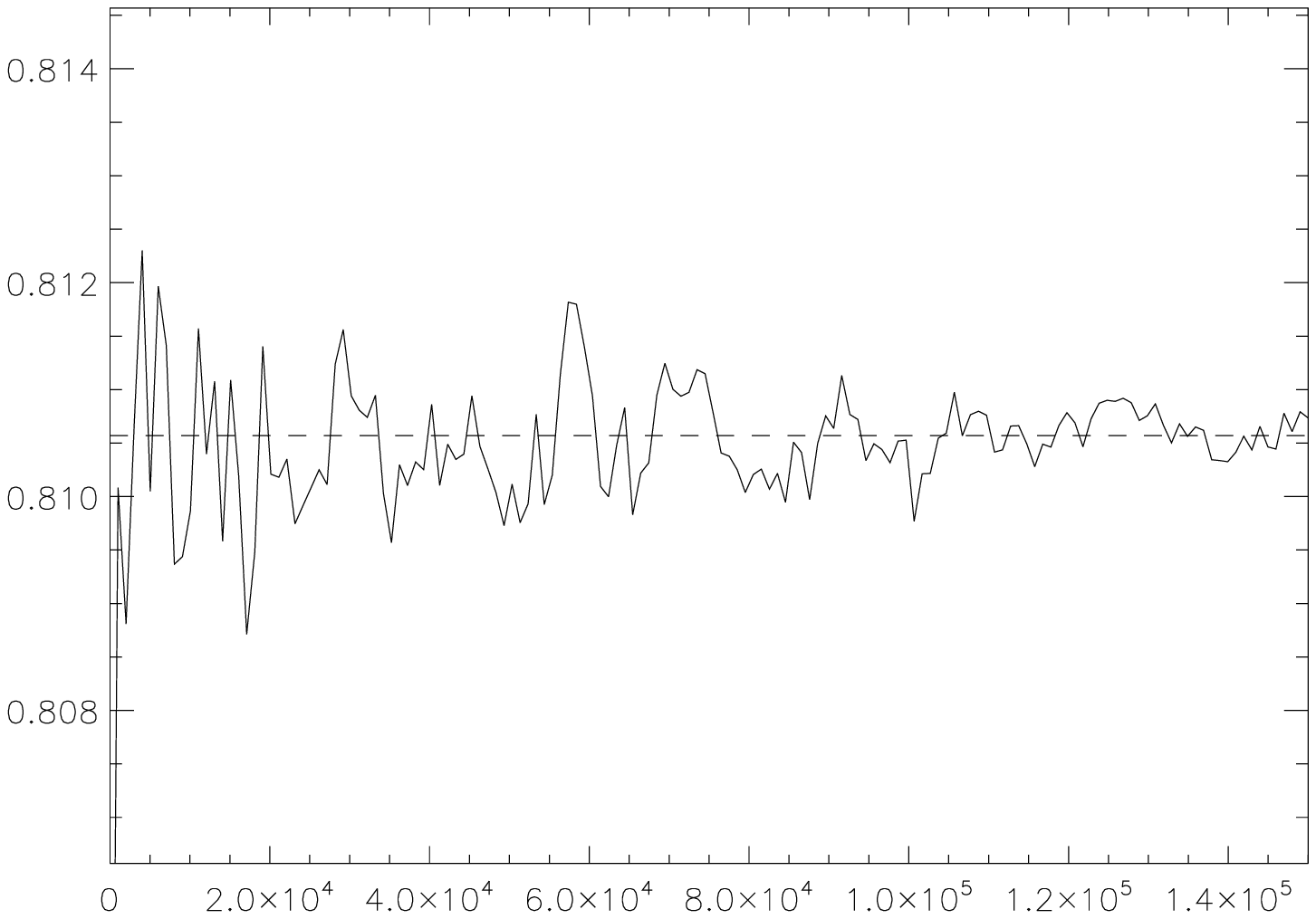}
\hspace{-0.5cm}\includegraphics[width=0.5\textwidth]{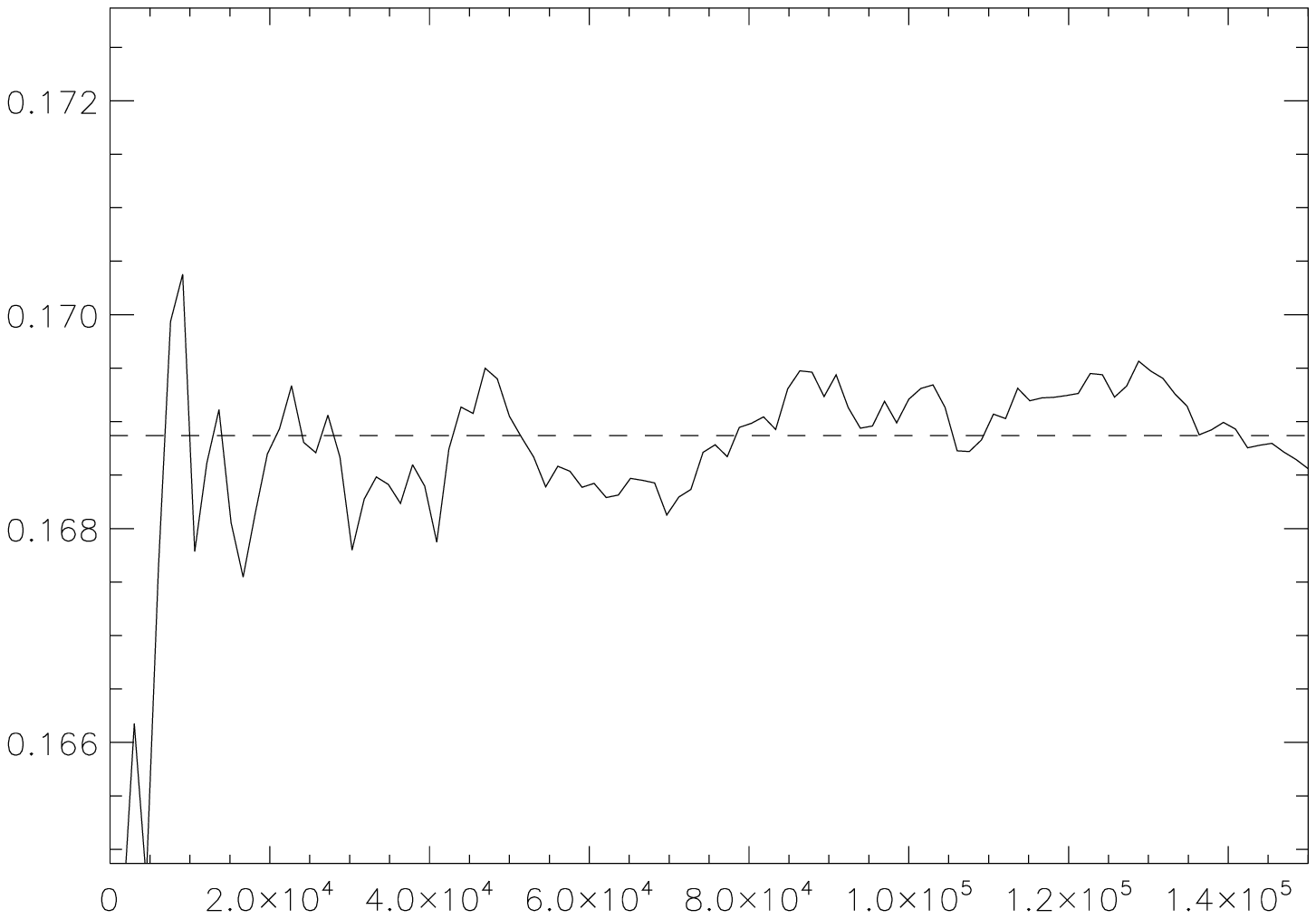}

\caption{{\it Left:} m=n=1, D=1, f=1. \ {\it Right:}  m=n=5, D=3, f=5.}
\label{fig:2}       
\end{figure}

\begin{figure}[!ht]
\hspace{-0.5cm}\includegraphics[width=0.5\textwidth]{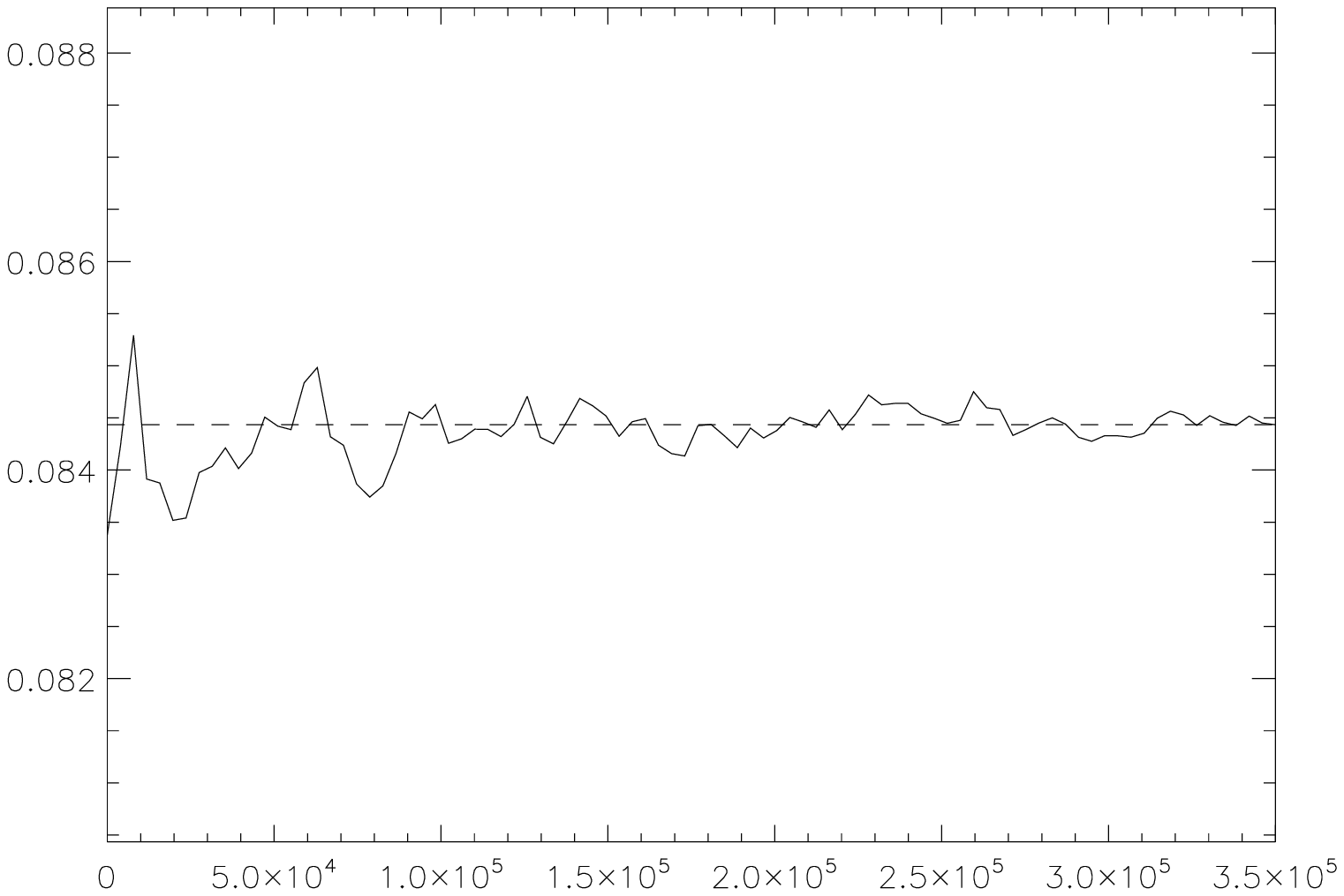}
\hspace{-0.5cm}\includegraphics[width=0.5\textwidth]{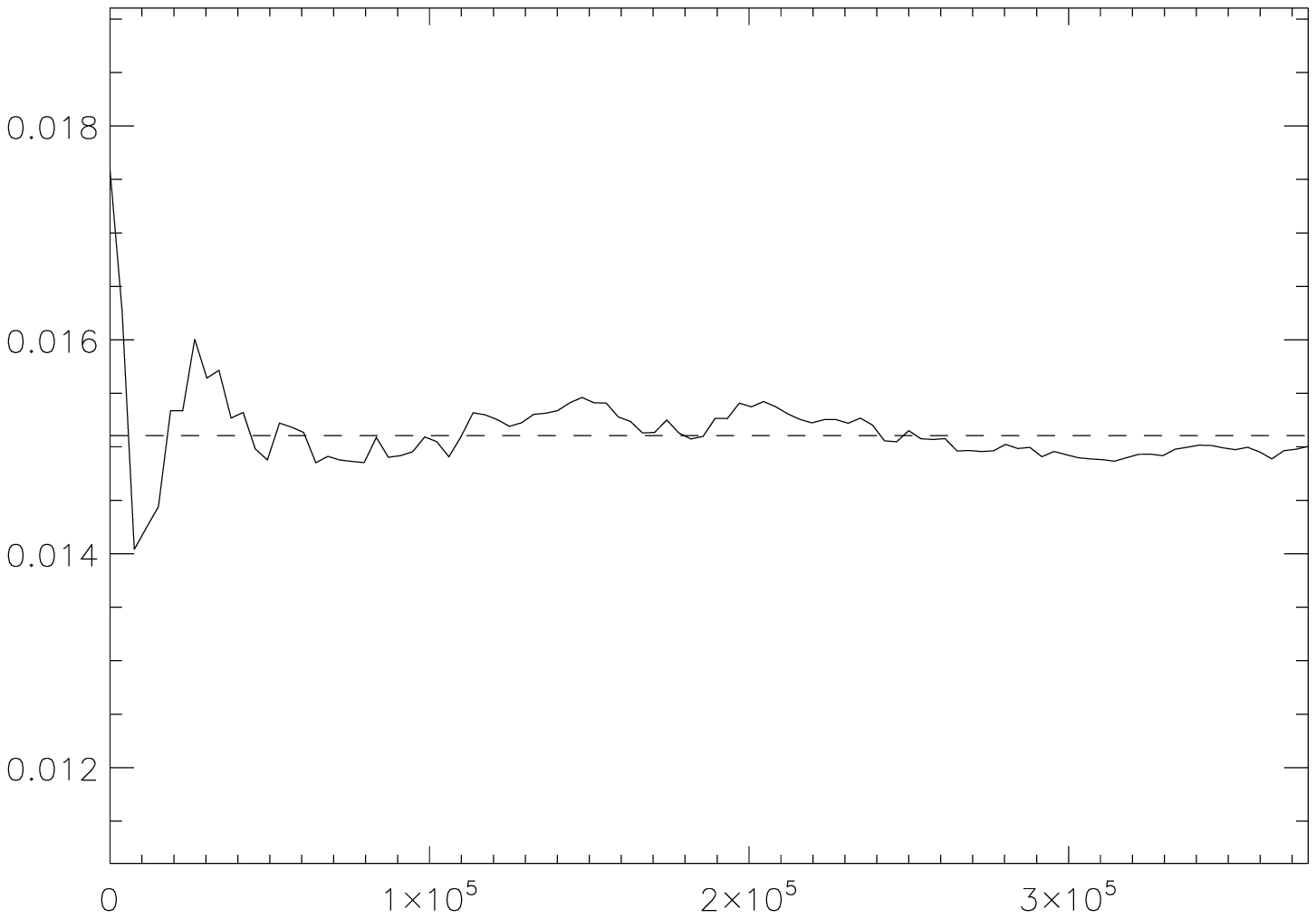}
\caption{{\it Left:} m=n=5, D=7, f=5. \ {\it Right:}  m=n=5, D=125, f=1.}
\label{fig:3}     
\end{figure}

For fixed $D$, $\chi$, $m$ and $n$, let 

\begin{equation*}
C(X)= \frac{1}{X} \sum_{\substack{0<c\leqslant X \\ c\congru 0\mmod D)}}
\frac{K_2(m,n;c)}{\sqrt{c}} \barre{\epsilon_c} \chi(c). 
\end{equation*}
The following examples (see Fig.\,\ref{fig:1}, Fig.\,\ref{fig:2} and Fig.\,\ref{fig:3}) illustrate the convergence of $C(X)$ to the value $C=C(D,f,m,n)$, given by Theorem\,\ref{intro:thm}, and shown in dotted line. Since the formula for $C(D,f,m,n)$ is quite involved, we consider several numerical examples. 

\begin{samepage}
Let us give more details about the exact value $C$ corresponding to each example and about the ranges for which $C(X)$ has ben computed:\\
\begin{itemize}

\item[-] Example\,1: $m=n=1$, $D=27$, $f=1$. Then $C= (3\pi^2)^{-1}=0.03377\pts$. \\
The function $C(X)$ is computed for $X\leqslant 2.7 \pt 10^5 $.

\item[-] Example\,2: $m=n=1$, $D=45$, $f=1$. Then $C= (6\pi^2)^{-1}=0.01688\pts$.\\
The function $C(X)$ is computed for $X\leqslant 2.7 \pt 10^5$. 

\item[-] Example\,3: $m=n=1$, $D=1$, $f=1$. Then $C= 8/\pi^2=0.81057\pts$. \\
The function $C(X)$ is computed for $X\leqslant 1.5 \pt 10^5$. 

\item[-] Example\,4:  $m=n=5$, $D=3$, $f=5$. Then $C= 5/(3\pi^2)=0.16886\pts$. \\
The function $C(X)$ is computed for $X\leqslant 1.5 \pt 10^5 $.

\item[-] Example\,5: $m=n=5$, $D=7$, $f=5$. Then $C= 5/(6\pi^2)=0.08443\pts$. \\
The function $C(X)$ is computed for $X\leqslant 3.5 \pt 10^5 $.

\item[-] Example\,6: $m=n=5$, $D=125$, $f=1$. Then $C= \sqrt{5}/(15\pi^2)=0.01510\pts$. 
The function $C(X)$ is computed for $X\leqslant 3.75 \pt 10^5 $.

\end{itemize}
\end{samepage}


\end{document}